\documentclass[12pt]{article}

\usepackage{amsfonts}
\usepackage{amsmath}
\usepackage{amssymb}
\usepackage{graphicx}

\newtheorem{theorem}{Theorem}[section]
\newtheorem{proposition}[theorem]{Proposition}
\newtheorem{corollary}[theorem]{Corollary}
\newtheorem{lemma}[theorem]{Lemma}
\newtheorem{remark}[theorem]{Remark}

 \def\bit{\begin{itemize}}
 \def\eit{\end{itemize}}

\def\oo{{\omega}}

\def\dd{\displaystyle}
\def\ooo{{\Omega}}
\def\<{\left<}
\def\>{\right>}
\def\({\left(}
\def\){\right)}
\def\ff{\forall }
\def\D{\Delta }
\def\9{{\infty}}
\def\barr{\begin{array}}
\def\earr{\end{array}}
\def\ov{\overline}
\def\wt{\widetilde}
\def\vp{{\varepsilon}}
\def\G{{\Gamma}}
\def\pp{{\partial}}
\def\dd{\displaystyle}
\def\vf{{\varphi}}
\def\lbb{{\lambda}}
\def\g{{\gamma}}
\def\a{{\alpha}}
\def\b{{\beta}}

\def\de{{\delta}}

\def\cald{\mathcal{D}}

\def\calh{\mathcal{H}}

\def\vsp{\vspace*{1,5mm}\\ }

\def\vspp{\vspace*{1,5mm}\\ }
\def\n{\noindent }
\def\eq{equa\-tion}

\def\inr{\mbox{ in }}
\def\onr{\mbox{ on }}
\def\D{{\Delta}}
\def\rr{{\mathbb{R}}}

\def\dd{\displaystyle}

\def\3{\subset }
\def\na{{\nabla}}

\def\mk{\medskip }
\def\bk{\bigskip }

\def\fwg{fol\-lo\-wing}

\title{\bf Nonlinear diffusion equations in~image processing}

\author{\bf Viorel Barbu}

\date{{\small
Al.I. Cuza University and Octav Mayer Institute of Mathematics\\ (Romanian Academy), Ia\c si, Romania}}

\begin{document}

\maketitle

\begin{abstract}
One surveys here a few nonlinear diffusion models in image restoration and denoising with main emphasis on that described by nonlinear parabolic equations of gradient type. The well-posedness of the corresponding Cauchy problem as well as stability of the derived finite difference scheme is studied from perspectives of nonlinear semigroup theory. As a matter of fact, most of denoising PDE procedures existing in literature though apparently are efficient at experimental level, are however mathematically ill posed and our effort here is to put them on more rigorous mathematical basis.
\end{abstract}

\section{Introduction}
Let us briefly review the PDE model for restoring the degraded (noisy) images in $\rr^d$. We consider the degraded image $u_0$ as a real valued function on a given bounded domain $\ooo$ of $\rr^d$ with the smooth boundary $\pp\ooo$. The restored (denoised) image is denoted by $u=u(x)$, $x\in\ooo$. The problem of recovering $u$ from $u_0$, via variational and regularization (smoothing) techniques, is a fundamental problem in image processing.

If $d=2$, then the domain $\ooo$ corresponds to two  dimensional images, while in the case $d=3$ we have three dimensional images which are of interest in medical imagistic. We note also that the case $u:\ooo\to\rr$  considered here corresponds to gray value images, while if $u:\ooo\to\rr^d$, $d=2,3$, the restoring process refers to colour images.

The classical denoising procedure is via linear heat equation (linear diffusion filter)
\begin{equation}\label{e1.1}
\barr{l}
\dd\frac{\pp u}{\pp t}=\D u\ \inr\ (0,\9)\times\rr^d,\vsp
u(0,x)=u_0(x),\ x\in\rr^d,
\earr\end{equation}
applicable also on bounded domains $\ooo$ if we take in \eqref{e1.1} the Dirichlet boundary condition $u=u^0$ on $\pp\ooo$ or the Neumann boundary condition $\frac{\pp u}{\pp\nu}=\na u\cdot\nu=0$ on $\pp\ooo$.

The solution to \eqref{e1.1} is given by
\begin{equation}\label{e1.2}
u(t,x)=(G(t)*u_0)(x)=\int_{\rr^d}G(t,x-y)u_0(y)dy,\end{equation}where $G$ is the Gaussian kernel
\begin{equation}\label{e1.3}
G(t,x)=\frac1{(2\sqrt{\pi t})^d}\ e^{-\frac{|x|^2}{4t}}.\end{equation}(Here,  $|\cdot|$  is the Euclidean norm of $\rr^d.$)

It should be noticed that this restoring method is related to the variational approach
\begin{equation}\label{e1.4}
{\rm Min}\left\{\int_\ooo|\na u|^2dx+\lbb\int_\ooo|u-u_0|^2dx,\ u\in H^1(\ooo)\right\},\end{equation}where $\lbb>0$ is a scale parameter.

Indeed, by the classical Dirichlet principle, the minimization problem \eqref{e1.4} considered on the Sobolev space $H^1(\ooo)$ is equivalent to the linear elliptic problem
\begin{equation}\label{e1.5}
-\D u+u=u_0\inr\ \ooo;\ \frac{\pp u}{\pp \nu}=0\ \onr\ \pp\ooo,\end{equation}while the solution $u$ to the parabolic boundary value problem \eqref{e1.1} with homogeneous Neumann's boundary conditions is given by$$u(t,x)=\lim_{h\to0}u_h(t,x),\ t\ge0,\ x\in\ooo,$$where $u_h$ is the solution to the finite difference scheme
\begin{equation}\label{e1.6}
\barr{l}
u_h^{i+1}(x)-h\D u^{i+1}_h(x)=u^i_h(x),\ i=0,1,...,\vsp
u^0_h=u_0,\ u_h(t)=u^i_h\mbox{ for }t\in(ih,(i+1)h).\earr\end{equation}In other words,
\begin{equation}\label{e1.7}
u^{i+1}_h={\rm arg}\min\left\{\int_\ooo|\na u|^2dx+\frac1{h}\int_\ooo|u-u^i_h|^2dx\right\}.\end{equation}The linear diffusion filter, though has a strong smoothing effect on the initial image $u_0$, has several disadvantages and  most important  is that it blurs the edges and so it makes hard their identification.
This limitation of linear filters was partially removed by nonlinear diffusion filters designed on similar variational principles as \eqref{e1.4}. Namely, the idea was to replace \eqref{e1.4} by
\begin{equation}\label{e1.8}
{\rm Min}\left\{\int_\ooo j(\na u)dx+\lbb\int_\ooo|u-u_0|^2dx\right\},\end{equation}where $j:\rr^d\to\rr$ is a given continuous function and $u$ is taken in an appropriate Sobolev space $W^{1,p}$ on $\ooo$ or, more generally, on the space $BV(\ooo)$ of functions with bounded variations on $\ooo$ (see Section 4 below).

Formally, the Euler--Lagrange conditions of optimality in \eqref{e1.8} are given by the \eq
\begin{equation}\label{e1.9}
\barr{ll}
-{\rm div}(\b(\na u))+2\lbb(u-u_0)=0& \inr\ \ooo,\vsp
(\na u)\cdot\nu=0& \onr\ \pp\ooo,\earr\end{equation}where $\b(r)=\na j(r),\ \ff r\in\rr.$

The corresponding evolution equation is
\begin{equation}\label{e1.10}
\barr{ll}
\dd\frac{\pp u}{\pp t}-{\rm div}(\b(\na u))=0&\inr\ (0,T)\times\ooo,\vsp
\b(\na u)\cdot\nu=0&\onr\ (0,T)\times\pp\ooo,\vsp
u(0)=u_0&\inr\ \ooo.\earr\end{equation}

If the minimum in formula \eqref{e1.8} is taken on the Sobolev space $W^{1,p}_0(\ooo)=\{u\in W^{1,p}(\ooo);\ u=0\ \onr\ \pp\ooo\}$, then in \eqref{e1.9} the Neumann boundary condition $\frac{\pp u}{\pp \nu}=\na u\cdot\nu=0
$ on $\pp\ooo$ is replaced by the homogeneous boundary condition $u=0$ on $\ooo$ with the corresponding modification in \eqref{e1.10}.

In specific applications to image restoring, the continuous image $u=u(x):\ooo\to\rr$ is viewed as a discrete image $u(ih,jh)=u_{ij}$, $i=1,...,N,$ $j=1,...,M$, where $u_{ij}\in\rr^+$ display the gray values at each pixel $(i,j)$ and $h$ is the $x$ grid spacing. Then the gradient $\na u$ at $x_i=ih$, $y_j=jh$ is approximated by

$$\{u_{i,j}-u_{i-1,j},u_{i,j}-u_{i,j-1}\}.$$(Here and everywhere in the following, by $\na u$ we mean the gradient of\break $x\to u(t,x).)$

It should be observed that \eqref{e1.8} implies, likewise \eqref{e1.9}, a smoothness constraint $j(\na u)\in L^1(\ooo)$, which is important for restoring  the initial degraded image $u_0$. In terms of the evolution equation \eqref{e1.10}, this property is implied by a smoothing effect on initial data $u_0$ of the flow  $t\to u(t)$ generated by differential equations of gradient types.

The first nonlinear filter of this form was proposed by Perona and Malik \cite{15}. Their model is of the form \eqref{e1.10}, where $$\barr{l}
j(r)=\dd\frac12\ j_0(|r|^2)\mbox{\ \ and}\vsp
j_0(s)=\a^2\log(s+\a^2),\ s\ge0.\earr$$That is,
\begin{equation}\label{e1.11}
\barr{ll}
\dd\frac{\pp u}{\pp t}-{\rm div}(g(|\na u|^2)\na u)=0&\inr\ (0,T)\times\ooo,\vsp
\na u\cdot\nu=0&\onr\ (0,T)\times\pp\ooo,\vsp
u(0)=u_0&\inr\ \ooo,\earr\end{equation}where
\begin{equation}\label{e1.12}
g(s)=\frac{\a^2}{\a^2+s},\ \ \ff s\ge0.\end{equation}The Perona-Malik model \eqref{e1.11}, \eqref{e1.12} was successful in many specific situations and it is at origin of a large variety of other models derived from \eqref{e1.11} (see, e.g., \cite{8aa}, \cite{10a}, \cite{10aa}, \cite{13}, \cite{17}, \cite{18} and the references given in these works).
It should be said, however, that a severe limitation of the model is that problem \eqref{e1.11} is, in general, ill posed. Indeed, by a simple computation it follows that the mapping $r\to g(|r|^2)r$ is monotonically decreasing in $\rr^d$ for $|r|>\a$  which, as is well-known, implies that the forward Cauchy problem \eqref{e1.11} is ill posed in this region. In order to avoid this inconvenience, which has as principal effect the instability of the finite difference scheme associated to \eqref{e1.11}, later on in literature,  equation \eqref{e1.11} was replaced by the regularized version (see, e.g., \cite{1}, \cite{17}, \cite{18})

\begin{equation}\label{e1.13}
\barr{ll}
\dd\frac{\pp u}{\pp t}-{\rm div}(g(|\na(G*u)|^2)\na u)=0,\ \inr\ (0,T)\times\ooo,\vsp
\na u\cdot\nu=0\ \onr\ (0,T)\times\pp\ooo,\ \ u(0)=u_0,\ \inr\ \ooo\earr\end{equation}where $G$ is the Gaussian kernel \eqref{e1.3}. A related version of \eqref{e1.13} is the  parabolic system
\begin{equation}\label{e1.14}
\barr{ll}
\dd\frac{\pp u}{\pp t}-{\rm div}(g(|\na v|^2)\na u)=0,&\inr\ (0,T)\times\ooo,\vsp
\dd\frac{\pp v}{\pp t}-\D v=u&\inr\ (0,T)\times\ooo,\vsp
\dd\frac{\pp u}{\pp\nu}=0,\ \dd\frac{\pp v}{\pp\nu}=0& \onr\ (0,T)\times\pp\ooo,\vsp
u(0)=u_0\ v(0)=0&\inr\ \ooo.\earr\end{equation}One can prove that this modified model given by \eqref{e1.13}, \eqref{e1.14} is well posed, but we omit  the details.
We shall see, however, in Section 3 that the original model \eqref{e1.11} is locally well posed in a sense to be made precise.  It should be said, also, that the Peron--Malik model \eqref{e1.11} as well as other restoring models of the form \eqref{e1.8}, which are defined on Sobolev spaces $W^{1,p}(\ooo)$, $1<p<\9$, are quite inefficient in edge detection and this fact has the \fwg\ simple mathematical explanation: a function $u\in W^{1,p}(\ooo)$ is continuous along rectifiable curves in $\ooo$ and so it cannot detect edges which represent discontinuities arising in image domain. For this reason, in \cite{6}  one chooses in \eqref{e1.8} $j:\rr^d\to\rr$ convex, continuous and such that
\begin{equation}\label{e1.15}
\lim_{|r|\to\9}\ \frac{j(r)}{|r|}=\lim_{|p|\to\9}\ \frac{j^*(p)}{|p|}=\9,\end{equation}where $j^*$ is the conjugate of $j$. In this way, the minimization problem \eqref{e1.8} reduces to $u\in W^{1,1}(\ooo)$, which provides in a certain sense a smooth restoring image with weaker continuity properties along rectifiable curves. From the point of view of edge detection, the best choice for \eqref{e1.8} (though hard to treat from the computational point of view) is {\it the total variation model} (BV model) introduced by L. Rudin et al. \cite{16} and developed in \cite{11},\cite{12} and \cite{14}. Namely,
\begin{equation}\label{e15}
{\rm Min}\left\{ \|Du\| +\lbb\int_\ooo|u-u_0|^2dx\right\},\end{equation}where $\|Du\|$ is a total variation of $u\in BV(\ooo).$

The advantage of the $BV$ model is that the functions with bounded varia\-tion might have discontinuity along rectifiable curves in $\ooo$ being, however, well defined
and even continuous on $\ooo$ in a generalized sense. This model will be discussed in Section 2.

\n{\bf Notation.} Let $\ooo$ be an open subset of $\rr^d$. Denote by $L^p(\ooo),\ 1\le p\le\9$, the space of Lebesgue measurable functions on $\ooo$ such that
$$|u|_p=\(\int_\ooo|u(x)|^p dx\)^{1/p}<\9,\ 1\le p<\9,$$with usual modification for $p=\9$. By $W^{1,p}(\ooo),\ 1\le p\le\9$, we denote the Sobolev space
$$W^{1,p}(\ooo)=\left\{u\in L^p(\ooo);\ \frac{\pp u}{\pp x_i}\in L^p(\ooo),\ i=1,...,d\right\},$$where $\frac\pp{\pp x_i}$ are taken in the sense of distributions on $\ooo$. The Sobolev space $W^{2,p}(\ooo)$ is similarly defined.  We denote by $W^{1,p}_0(\ooo)$ the space $\{u\in W^{1,p}(\ooo);$ $ \g(u)=0$ on $\pp\ooo\}$, where $\g(u)$ is the trace of $u$ to $\pp\ooo$ -- the boundary of $\ooo$.  Everywhere in the following $\pp\ooo$ is assumed to be piecewise smooth in order to have a $W^{2,2}(\ooo)$ regularity for solutions to the elliptic boundary value problem.   We set $H^1(\ooo)=W^{1,2}(\ooo),$ $H^1_0(\ooo)=W^{1,2}_0(\ooo).$

Given a Banach space $X$ with the norm $\|\cdot\|_X$, we denote by $C([0,T];X)$ the space of all continuous $X$-valued functions $u:[0,T]\to X$. By $L^p(0,T;X)$, $1\le p\le\9$, we denote the space of all Bochner measurable $L^p$-integrable functions $u:(0,T)\to X$ with the norm
$$\|u\|_{L^p(0,T;X)}=\(\int^T_0\|u(t)\|^p_Xdt\)^{1/p}.$$By $W^{1,p}([0,T];X)$, we denote the space of all absolutely continuous functions $u:[0,T]\to X$ which are a.e. differentiable on $(0,T)$ and $\frac{du}{dt}\in L^p(0,T;X)$. We refer to \cite{8}, \cite{9a} and \cite{9} for other notations and basic results on PDE relevant in the \fwg.

\section{Review on nonlinear semigroups\\ and~evolution equations in Hilbert spaces}
\setcounter{theorem}{0}
\setcounter{equation}{0}

Let $H$ be a real Hilbert space with the scalar product $\<\cdot,\cdot\>$ and norm $|\cdot|_H$. The nonlinear multivalued operator $A:D(A)\subset H\to H$ is said to be {\it monotone} if
$$\<y_1-y_2,x_1-x_2\>\ge0,\ \ff y_i\in Ax_i,\ i=1,2.$$

The operator $A$ is said to be {\it maximal monotone} if it is monotone and \mbox{$R(I+A)=H$.} Here, $R$ is the range and $D(A)=\{x\in H;\ Ax\ne\emptyset\}.$

\n Given a convex and lower semicontinuous function $\vf:H\to\ov\rr=\mbox{$]-\9,+\9]$}$, denote by $\pp\vf:H\to H$ the {\it subdifferential} of $\vf$, that is,
$$\pp\vf(x)=\{y\in H;\ \vf(x)\le\vf(u)+\<y,x-u\>,\ \ff u\in H\}.$$The operator $\pp\vf$ is maximal monotone (see, e.g., \cite{8}, \cite{9}).

Given $u_0\in H$ and $f\in L^1(0,T;H),\ 0<T<\9$, we consider the Cauchy problem
\begin{equation}\label{e2.1}
\barr{ll}
\dd\frac{du}{dt}\,(t)+Au(t)\ni f(t),\ t\in(0,T),\vsp
u(0)=u_0.\earr\end{equation}

By strong solution to \eqref{e2.1} we mean an absolutely continuous function $u:[0,T]\to H$ such that
\begin{equation}\label{e2.2}
\barr{ll}
\dd\frac{du}{dt}\,(t)+w(t)=f(t),&\mbox{a.e.}\ t\in(0,T),\vsp
w(t)\in Au(t),&\mbox{a.e.}\ t\in(0,T),\vsp
u(0)=u_0.\earr\end{equation}The main existence result for the Cauchy problem \eqref{e2.1}, which is valid in a general reflexive Banach space $H$ is the \fwg\ (see \cite{8},  p.~127, and~\cite{9}).

\begin{theorem}\label{t2.1} Let $A$ be maximal monotone and $f\in W^{1,1}([0,T];H)$, $u_0\in D(A)$. Then, the Cauchy problem \eqref{e2.1} has a unique strong solution $u=u(t,u_0,f)\in W^{1,\9}([0,T;H).$

Moreover, $y$ is everywhere differentiable from the right on $[0,T]$ and
$$\frac{d^+}{dt}\,y(t)+(Ay(t)-f(t))^0=0,\ \ff t\in[0,T),$$where $(Ay(t)-f(t))^0$ is the minimal section of the set $\{z;\ z\in Ay(t)-f(t)\}.$

The mapping $(u_0,f)\to u(\cdot,u_0,f)$ is Lipschitz from $H\times L^1(0,T;H)$ to $C([0,T];H)$, that is,
\begin{equation}\label{e2.3a}
|u(\cdot,u_0,f)-u(\cdot,\bar u_0,\bar f)|_H\le
|u_0-\bar u_0|_H+\dd\int^T_0|f(s)-\bar f(s)|_Hds,\end{equation}for all $u_0,\bar u_0\in D(A),\ f,\bar f\in W^{1,1}([0,T];H)$.
\end{theorem}

By \eqref{e2.3}, we see that $u=u(t,u_0,f)$ can be extended by continuity to all $u_0\in \ov{D(A)}$ and $f\in L^1(0,T;H)$. This extension is called a generalized solution to problem \eqref{e2.1} and, in general, it is not a strong solution, that is, a.e. differentiable on $(0,T).$

In the special case $f\equiv0$, \eqref{e2.1} defines a semiflow (semigroup) $S(t)u_0=u(t)$, $t\ge0$, which is a semigroup of contractions on $\ov{D(A)}$, that is,
$$\barr{ll}
S(t)S(s)u_0=S(t+s)u_0,&\ff t,s\ge0,\vsp
|S(t)u_0-S(t)\bar u_0|_H\le|u_0-\bar u_0|_H,&\ff u_0,\bar u_0\in\ov{D(A)},\vsp
\dd\lim_{t\to0} S(t)u_0=S(0)u_0=u_0,&\ff u_0\in\ov{D(A)}.\earr$$Moreover, $S(t)u_0$ is given by the exponential formula
\begin{equation}\label{e2.3}
S(t)u_0=\lim_{n\to\9}\(I+\frac tn\,A\)^{-n} u_0,\ \ff t\ge0.\end{equation}This means that the solution $u$ to \eqref{e2.1} (with $f\equiv0$) is given as the limit of finite difference scheme
\begin{equation}\label{e2.4}
\barr{ll}
u^h(t)=u^h\mbox{ \ \ for }t\in[ih,(i+1)h),\vsp
u^h_i+hAu^h_i\in u^h_{i-1},\ \ i=1,2,...,N=\left[\frac Th\right],\\
u^h_0=u_0.\earr\end{equation}It turns out that, in the special case where $A$ is a subgradient operator, that is, $A=\pp\vf$, where $\vf:H\to\ov\rr$ is a lower semicontinuous, convex function, the nonlinear semigroup $S(t)$ has a smoothing effect on initial data. Namely, one has (see \cite{8}, p.~170, \cite{9}).

\begin{theorem}\label{t2.2} Let $A=\pp\vf$ and $u_0\in\ov{D(A)}$. Then, the generalized solution $u(t)=S(t)u_0$ to \eqref{e2.1} is a strong solution and $u\in W^{1,\9}([\de,T);H)$, for all $0<\de<T$,
\begin{equation}\label{e2.4a}
\sqrt{t}Au\in L^2(0,T;H),\ t\vf(u(t))\in C([0,T]).\end{equation}In particular, it follows that
\begin{equation}\label{e2.5}
\barr{ll}
\dd\frac{d^+u}{dt}\,(t)+(\pp\vf(u(t)))^0=0,\ \ \ff t\in(0,T),\vsp
u(0)=u_0.\earr\end{equation}\end{theorem}
Theorem \ref{t2.2} extends to the nonhomogeneous equation \eqref{e2.1} with $f\in L^2(0,T;H)$.\newpage

By \eqref{e2.4}, it follows that, for each $u_0\in\ov{D(A)}$, $$S(t)u_0=u(t,u_0)\in D(A),\ \ff t>0,$$ and this means that the semigroup $S(t)$, also denoted $e^{tA}$, has a smoothing effect on initial  data. Since, in applications to partial differential equations, $H$ is usually the space $L^2(\ooo)$, while $D(A)$ is a
Sobolev space on $\ooo$ (very often $H^2(\ooo))$, which incorporates certain boundary conditions, the significance of "smoothing" effect becomes apparent.

In applications to image restoring,   $u_0$ is the degraded image, usually in $L^2(\ooo)$, while $u=u(t,x)$, $x\in\ooo$, is the recovered image.

The fact that $x\to u(t,x)$ is smooth, though $u(0)=u_0(x)$ is not, is of crucial importance in restoring (denoising) processes. It should be emphasized, however, that, from all maximal monotone nonlinear operators $A$ on $H$, only subgradients mappings have this smoothing property.

We note also that this image restoring process is related and, in a certain sense,  equivalent with the variational method
\begin{equation}\label{e2.6}
{\rm Min}\{\vf(u)+\lbb|u-u_0|^2_H\}.\end{equation}Equivalently,
\begin{equation}\label{e2.7}
\pp\vf(u)+2\lbb(u-u_0)=0.\end{equation}Indeed, as seen above, the solution $u=u(t)$ is given by $u(t)=\dd\lim_{h\to0}u_h(t)$, where $u_h$ is the solution to the difference scheme \eqref{e2.4}, that is,
$$u^h_i+hAu^h_i\ni u^h_{i-1},\ i=1,...$$Equivalently,
\begin{equation}\label{e2.8}
u^h_i={\rm arg}\,{\rm Min}\left\{\vf(v)+\frac1{2h}\,|v-u^h_{i-1}|^2_H\right\},\ i=1,...\end{equation}It should be said that both approaches, \eqref{e2.1} and \eqref{e2.6}, have a smoothing effect on the initial image $u_0$. However, we see by \eqref{e2.8} that the evolution process has a cumulative smoothing effect via the iterative process $u^h_i$.  In the linear case, for instance if $A$ is a linear elliptic operator, we see that $u(t)\in\dd\bigcap_{k\ge1}D(A^k)=C^\9(\ov\ooo)$, which, as seen earlier, precludes the identification of edges. For nonlinear operators, in  general, this does not happen, the maximal regularity of $u(t)$ being that given by $u(t)\in D(A)$.  This fact explains the advantage of nonlinear filters of the form \eqref{e2.5}, which keep a certain equilibrium between the   image smoothing and the edge detection.

For the denoising procedure, it suffices to consider equation \eqref{e2.5} or its discrete version \eqref{e2.8} on an arbitrary small interval $[0,T]$. As a matter of fact, considering equation \eqref{e2.5} on a large interval $[0,T]$ might alter the result  because, for $t\to\9$, the solution $u=u(t)$ to \eqref{e2.5} is weakly convergent to an equilibrium solution $u_\9$ to \eqref{e2.5}, that is, $\pp\vf(u_0)\ni0$. (See \cite{8}, p.~172.) Hence, the best results are obtained for $T$ small.

However, by virtue of the same  asymptotic result, we can recover the blurred image $u_0$ from the evolution equation
\begin{equation}\label{e2.9}
\barr{ll}
\dd\frac{du}{dt}\,(t)+Au(t)+2\lbb(u(t)-u_0)\ni0,&\ff t\ge0,\vsp
u(0)=u^0,\earr\end{equation}where $A=\pp\vf$, $u^0\in D(A)$. Indeed, as mentioned above, for $t\to\9$, the solution $u(t)$ to \eqref{e2.9} is weakly convergent in $H$ to a solution $u_\9$ to the stationary equation
$$Au_\9+2\lbb(u_\9-u_0)\ni0,$$which is just a minimum point for problem  \eqref{e2.6}. The dynamic model \eqref{e2.9} is equivalent with the classical steepest descent algorithm for the minimization problem \eqref{e2.6}. Namely,
$$u_{i+1}+h(A(u_{i+1})+2\lbb(u_{i+1}-u_0))=u_i.$$

\section{Restoring the image via nonlinear parabolic equations in divergent form}
\setcounter{theorem}{0}
\setcounter{equation}{0}

We consider here the restoring model given by the parabolic boundary value problem
\begin{equation}\label{e3.1}
\barr{ll}
\dd\frac{\pp u}{\pp t}-{\rm div}(\b(\na u))=0&\inr\ (0,T)\times\ooo,\vsp
\b(\na u)\cdot\nu=0&\onr\ (0,T)\times\pp\ooo,\vsp
u(0,x)=u_0(x),&x\in\ooo.\earr\end{equation}Here, $\b=\na j$, where $j:\rr^d\to\rr$, $d=2,3$, is a convex and differentiable function satisfying the growth conditions
\begin{equation}\label{e3.2}
\oo_1|r|^{p}+C_1\le j(r)\le\oo_2|r|^p+C_2,\ \ff r\in\rr^d,\end{equation}for some $\oo_i>0,\ C_i\in\rr,\ i=1,2.$ Here, $1<p<\9.$

It is easily seen that $\b:\rr^d\to\rr^d$ is monotone, continuous and
$$\barr{rcll}
|\b(r)|&\le&C_3(|r|^{p-1}+1),&\ff r\in\rr,\vsp
\b(r)\cdot r&\ge&\oo_1|r|^p+C_4,&\ff r\in\rr,\earr$$where $C_3,C_4\in\rr.$

Consider in the space $H=L^2(\ooo)$ the nonlinear operator
\begin{equation}\label{e3.3}
\!\!\barr{l}
Au=-{\rm div}(\b(\na u)),\ \ \ff u\in D(A)\subset H,\vsp
D(A)=\{u\in W^{1,p}(\ooo);\, {\rm div}\,\b(\na u)\in L^2(\ooo),\, \b(\na u)\cdot\nu=0\ \onr  \pp\ooo\}.\earr\end{equation}(Here, $\b\cdot(\na u)\cdot\nu=0$ on $\pp\ooo$ should be understood in the \fwg\ weak sense
$$\int_\ooo\b(\na u)\cdot\na\psi\,dx=-\int_\ooo{\rm div}\,\b(\na u)\psi\,dx,\ \ff\psi\in C^1(\ov\ooo).)$$We have

\begin{theorem}\label{t3.1} The operator $A$ is maximal monotone in $L^2(\ooo)$. More precisely, $A=\pp\vf$, where $\vf:L^2\ooo)\Longrightarrow \ov\rr$,
\begin{equation}\label{e3.4}
\vf(u)=\left\{\barr{ll}
\dd\int_\ooo j(\na u)dx&\mbox{if }u\in W^{1,p}(\ooo),\ j(\na u)\in L^1(\ooo),\vsp
+\9&\mbox{otherwise}.\earr\right.\end{equation}\end{theorem}

\n{\bf Proof.} It is easily seen that $\vf:L^2(\ooo)\to\ov\rr$ is convex and lower-semicontinuous. Moreover, if $\eta\in Au$, we have
$$\barr{lcl}
\dd\int_\ooo Au(u-v)dx&=&\dd\int_\ooo\b(\na u)\cdot(\na u-\na v)dx\vsp
&\ge&\dd\int_\ooo(j(\na u)-j(\na v))dx,\ \ff v\in L^2(\ooo),\earr$$which implies that $A\subset\pp\vf$. In order to prove that $A=\pp\vf$, it suffices to show that $A$ is maximal monotone, that is, $\rr(I+A)=L^2(\ooo)$. To this end, we consider, for $f\in L^2(\ooo)$, the elliptic  equation
\begin{equation}\label{e3.5}
\barr{ll}
u-{\rm div}\, \b(\na u)=f&\inr\ \ooo,\vsp
\b(\na u)\cdot\nu=0&\onr\ \pp\ooo.\earr
\end{equation}(Here and everywhere in the following, div and $\na$ are taken in the sense of distributions on $\ooo$).

We associate with \eqref{e3.5} the minimization problem
\begin{equation}\label{e3.6}
{\rm Min}\left\{\Phi(u)=\dd\int_\ooo j(\na u)dx+\dd\frac12\int_\ooo(u-f)^2dx;\ u\in W^{1,p}(\ooo)\right\}.\end{equation}The function $\Phi$ is convex, lower semicontinuous in $L^2(\ooo)$  and by \eqref{e3.2} we see that the level sets $\{u;\ \Phi(u)\le\lbb\}$ are bounded in $W^{1,p}(\ooo)\cap L^2(\ooo)$ and, therefore, weakly compact. This implies that $\Phi$ attains its infimum at $u^*\in W^{1,p}(\ooo)\cap L^2(\ooo)$. We have, therefore,
$$\Phi(u^*)\le\Phi(u^*+\lbb v),\ \ff v\in C^1(\ov\ooo)$$and this yields
$$\int_\ooo\b(\na u^*)\cdot\na v\,dx+\int_\ooo(u^*-f)v\,dx=0,\ \ff v\in C^1(\ov\ooo).$$Hence, $u^*$ is the solutions to \eqref{e3.5}, as desired.

Then, by Theorem \ref{t2.2}, we obtain

\begin{theorem}\label{t3.2} For each $u_0\in L^2(\ooo)$, the equation
\begin{equation}\label{e3.7}
\barr{ll}
\dd\frac{\pp u}{\pp t}-{\rm div}\,\b(\na u)=f&\inr\ (0,T)\times\ooo,\vsp
\b(\na u)\cdot\nu=0&\onr\ (0,T)\times\pp\ooo,\vsp
u(0)=u_0,&\inr\ \ooo,\earr\end{equation}has a unique solution $u\in C([0,T];L^2(\ooo))$ such that
\begin{eqnarray}\label{e3.8}
\sqrt{t}\,\dd\frac{\pp u}{\pp t}\,,\sqrt{t}\ {\rm div}\,\b(\na u)\in L^\9(0,T;L^2(\ooo)),\\[2mm]
t\int_\ooo|\na u(t,x)|^pdx\le C,\ \ff t\in[0,T].\label{e3.9}
\end{eqnarray}\end{theorem}

By \eqref{e3.8}, \eqref{e3.9}, it follows the smoothing effect on the original image $u_0$.

Theorem \ref{t3.2} remains true in the anisotropic case $\b=\b(x,\cdot)$, where $\b(x,r)=\na j(x,r)$ and $j:\ooo\times\rr\to\rr$ satisfies conditions \eqref{e3.2} uniformly with  respect to $x\in\ooo$.

The case $p=1$ was excluded here. However, if $j$ satisfies conditions \eqref{e1.15}, one obtains a similar result with $u(t,\cdot)\in W^{1,1}(\ooo)$ (see \cite{8}, p. 81, \cite{6}). As mentioned earlier from point of view of edges detection smaller $p$ are more convenient in problem \eqref{e3.6} or \eqref{e3.7}.

\section{The total variation flow approach to image restoring}
\setcounter{theorem}{0}
\setcounter{equation}{0}

The literature on the total variation flow method in image processing contains a very large number of paper following to the pioneering work of L. Rudin, S. Osher and E. Fatemi \cite{16}. Among the most important contributions in this field, the works \cite{11}, \cite{12}, \cite{14} must be cited in first lines.

Formally, the evolution equation  defining the total variation flow is given~by \begin{equation}\label{e4.1}
\barr{ll}
\dd\frac{\pp u}{\pp t}-{\rm div}\(\frac{\na u}{|\na u|}\)=0&\inr\ (0,T)\times\ooo,\vsp
\na u\cdot\nu=0&\onr\ (0,T)\times\pp\ooo,\vsp
u(0)=u_0&\inr\ \ooo,\earr\end{equation}or its stationary counterpart
\begin{equation}\label{e4.2}
{\rm Min}\left\{\int_\ooo|\na u(x)|dx+\lbb\int_\ooo|u(x)-u_0(x)|^2dx;\ u\in W^{1,1}(\ooo)\right\}.\end{equation}

However, the parabolic boundary value problem \eqref{e4.1} as well as the mi\-ni\-mi\-za\-tion problem \eqref{e4.2} are ill posed because the nonlinear operator\break $u\to{\rm div}\(\frac{\na u}{|\na u|}\)$  with Neumann boundary condition (or with Dirichlet boun\-dary conditions as well) is not maximal monotone in $L^2(\ooo)$, while the energy functional $\Phi_0:L^2(\ooo)\to\ov\rr,$
$$\Phi_0(u)=\int_\ooo|\na u|dx+\lbb\int_\ooo|u(x)-u_0(x)|^2dx,\ \ff u\in W^{1,1}(\ooo),$$is not lower semicontinuous. Hence, $\Phi_0$ does not attains its infimum on $W^{1,1}(\ooo)$, but in a larger space, namely $BV(\ooo)$.

Let us now briefly review the definitions and basic properties of functions of bounded variation on $\ooo$. (For a general presentation of this space, we refer to \cite{2}, \cite{5a}.)

A function $u\in L^1(\ooo)$ is called a {\it function of bounded variations} if
\begin{equation}\label{e4.3}
\|Du\|=\sup\left\{\int_\ooo u\,{\rm div}\,\psi\,dx;\ \psi\in C^\9_0(\ooo;\rr^d),\,|\psi|_\9\le1\right\}<+\9.\end{equation}Equivalently, $\frac{\pp u}{\pp x_i}\in M(\ooo),\ i=1,...,d$, where $M(\ooo)$ is the space of Radon measures on $\ooo$. The space of functions of bounded variations on $\ooo$ is denoted by $BV(\ooo)$ and it is a Banach space with the norm

$$\|u\|_{BV(\ooo)}=|u|_1+\|Du\|.$$ $\|Du\|$ is called the {\it total variation} of $u\in BV(\ooo)$ and is also denoted by
$$\|Du\|=\int_\ooo|Du|,$$where $Du=\na u$ is the gradient of $u$ in sense of distributions.

For $u\in W^{1,1}(\ooo)$, we have $\|Du\|=\int_\ooo|\na u|dx$.

It turns out that $u\to \|Du\|$ is lower semicontinuous in $L^1(\ooo)$, that is,
\begin{equation}\label{e4.4}
\liminf_{u_n\buildrel{L^1(\ooo)}\over\longrightarrow u}\|Du_n\|\ge\|Du\|.\end{equation}Moreover, $BV(\ooo)\subset L^p(\ooo)$ for $1\le p\le\frac d{d-1}$ and, if $1\le p<\frac d{d-1}$, the above embedding is compact.

If $\pp\ooo$ is smooth (Lipschitz, for instance), there is a linear continuous operator $\g_0:BV(\ooo)\to L^1(\pp\ooo;\calh^{d-1})$ such that
$$\int_\ooo u\,{\rm div}\,\psi\,dx=-\int_\ooo\psi \cdot Du+\int_{\pp\ooo}(\psi\cdot\nu)\g_0(u)d\calh^{d-1},\ \ff\psi\in C^1(\ov\ooo;\rr^d).$$
(Here, $\calh^{d-1}$ is the Hausdorff measure of dimension $d-1$.) The function $\g_0(u)$ is called the {\it trace} of $u$ to $\pp\ooo$.

Consider the function $\vf:L^2(\ooo)\to\rr^d$ defined by
\begin{equation}\label{e4.5}
\vf(u)=\left\{\barr{ll}
\dd\int_\ooo|Du|&\mbox{if }u\in BV(\ooo)\cap L^2(\ooo),\vsp
+\9&\mbox{otherwise}.\earr\right.\end{equation}As easily seen, $\vf$ is convex and, by \eqref{e4.4}, it is lower semicontinuous on $L^2(\ooo)$. It turns out (see, e.g, \cite{5a}, p.~437) that the function $\vf$ is the closure $\ov\vf_0$ in $L^1(\ooo)$ of the function$\vf_0:L^1(\ooo)\to\ov\rr,$
\begin{equation}\label{e4.6}
\vf_0(u)=\left\{\barr{ll}
\dd\int_\ooo|\na u|dx&\mbox{if }u\in W^{1,1}(\ooo),\vsp
+\9&\mbox{otherwise},\earr\right.\end{equation}that  is, if $u_n\to u$ in $L^1(\ooo)$, then $$\vf(u)\le\dd\liminf_{n\to\9}\vf_0(u)$$ and there is $\{u_n\}\subset L^1(\ooo)$ convergent to $u$ in $L^1(\ooo)$, such that $$\dd\limsup_{n\to\9}\vf_0(u_n)\le\vf(u).$$This fact shows that $\vf$ is the natural lower semicontinuous extension of the energy function $\vf_0$ to $BV(\ooo)$. Moreover, since $\vf$ is convex and lower semicontinuous on $L^2(\ooo)$, the minimization problem
\begin{equation}\label{e4.7}
{\rm Min}\left\{\int_\ooo|Du|+\lbb\int_\ooo|u-u_0|^2dx;\ u\in L^2(\ooo)\right\}\end{equation}
has, for each $u_0\in L^2(\ooo)$ and $\lbb>0$, a unique solution $u^*\in BV(\ooo)\cap L^2(\ooo)$. Moreover, $u^*$ is the solution to the equation
\begin{equation}\label{e4.8}
\pp\vf(u^*)+2\lbb(u^*-u_0)\ni0,\end{equation}where $\pp\vf:L^2(\ooo)\to L^2(\ooo)$ is the subdifferential of $\vf$, that is,
\begin{equation}\label{e4.9}
\barr{r}
\pp\vf(u^*)=\{\eta\in L^2(\ooo);\ \vf(u^*)\le\vf(v)+\<\eta,u^*-v\>,\vsp \ff v\in BV(\ooo)\cap L^2(\ooo)\}.\earr\end{equation}We note that, if $u\in W^{1,1}(\ooo)$ and ${\rm div}\(\frac{\na u}{|\na u|}\)\in L^1(\ooo)$, $\na u\cdot\nu=0$ on $\pp\ooo$, then $u\in D(\pp\vf)$ and $\pp\vf(u)=-{\rm div}\(\frac{\na u}{|\na u|}\)$.

In general, the structure of $\pp\vf$ is more complex. One has, however, the \fwg\ description of the subdifferential $\pp\vf$ (see \cite{3}, \cite{4}, \cite{5}): $\eta\in \pp\vf(u)$ if and only if $u\in L^2(\ooo)\cap BV(\ooo)$, and there is $z\in L^\9(\ooo;\rr^d)$, ${\rm div}\,z\in L^2(\ooo)$, such that $\eta=-{\rm div}\,z$ and
\begin{equation}\label{e4.10}
\int_\ooo u\eta\,dx=\int_\ooo|Du|,\ \int_{\pp\ooo}(z\cdot\nu)\g_0(u)d\calh^{d-1}=0.\end{equation}

Consider now the evolution associated to $\pp\vf$,
\begin{equation}\label{e4.11}
\barr{l}
\dd\frac{d}{dt}\,u(t)+\pp\vf(u(t))\ni0,\ \ \ff t\in[0,T],\vsp
u(0)=u_0.\earr\end{equation}By virtue of the above discussion, the Cauchy problem \eqref{e4.11} is the "minimal" extension to problem \eqref{e4.1} to the space $BV(\ooo)$, where it is well posed. It should be said that in \eqref{e4.11}, likewise in the stationary problem \eqref{e4.8}, the Neumann boundary condition $\na u\cdot\nu=0$ are implicitly  incorporated in the condition $u(t)\in D(\pp\vf)$  through definition \eqref{e4.5} of $\vf$, the description \eqref{e4.10} of $\pp\vf$, and by the fact that $\vf=\ov\vf_0$.

By Theorem \ref{t2.1} and Theorem \ref{t2.2}, we have

\begin{theorem}\label{t4.1} For each $u_0\in L^2(\ooo)=\ov{D(\pp\vf)}$ there is a unique function $u\in C([0,T];L^2(\ooo))\cap W^{1,\9}([\de,T];H),$ $ \ff\de\in(0,T)$ such that
\begin{equation}\label{e4.12}
\barr{l}\dd\frac{d^+}{dt}\,u(t)+(\pp\vf(u(t)))^0=0,\ \ff t\in(0,T),\vsp
u(0)=u_0,\earr \end{equation}where $\pp\vf(u))^0$ is the minimal section of $\pp\vf(u)$. If $u_0\in D(\pp\vf)$, then $u\in W^{1,\9}([0,T];H)$.\end{theorem}

The function $u(t)=e^{t\pp\vf}u_0$ is the {\it total variation flow  on $\ooo$} and by \eqref{e4.12} it follows the smoothing effect of $e^{t\pp\psi}$ on initial data, that is, $u(t)\in BV(\ooo)$ for all $t\ge0$.

We may rewrite \eqref{e4.12} in terms of \eqref{e4.10}, but an explicit form of equation \eqref{e4.12} is, however, hard to find.

As seen earlier, the solution $u$ to \eqref{e4.12} can be obtained by the finite difference implicit scheme
$$u(t)=\lim_{h\to0}u_h(t)\mbox{ uniformly on $[0,T]$},$$where
\begin{equation}\label{e4.13}
\barr{ll}
u_h(t)=u^i_h\ \ \mbox{for }t\in[ih,(i+1)h),\ Nh=T,\vsp
u^{i+1}_h+h\pp\vf(u^{i+1}_h)\ni u^i_h,\ \ i=0,1,...,N.\earr\end{equation}However, solving directly \eqref{e4.13} is really difficult because of the complicate structure of $\pp\vf$. So, it is convenient to replace $\vf$ (or $\pp\vf)$ by more regular functions.
One such approximation of $\vf$ is described below. Namely, $\vf_\vp:L^2(\ooo)\to\ov\rr$ is defined by
\begin{equation}\label{e4.14}
\vf_\vp(u)=\left\{\barr{ll}
\dd\int_\ooo\(j_\vp(|\na u|)+\dd\frac\vp2\,|\na u|^2\)dx&\mbox{for }u\in H^1(\ooo),\vsp
+\9&\mbox{for }u\in L^2(\ooo)\setminus H^1(\ooo).\earr\right.\end{equation}Here, $j_\vp :[0,\9)\to[0,\9)$ is defined by
$$j_\vp(u)=\left\{\barr{ll}
\dd\frac1{2\vp}\,r^2&\mbox{for }0\le r\le\vp,\vsp
r-\dd\frac\vp2&\mbox{for }r>\vp.\earr\right.$$
Then the subdifferential $\pp\vf_\vp$ is given by
\begin{equation}\label{e4.15}
\barr{lcll}
\pp\vf_\vp(u)&=&-{\rm div}(\psi_\vp(\na u))-\vp\D u,\ \ff u\in D(\pp\vf_\vp),\vsp
D(\pp\vf_\vp)&=&\{u\in H^2(\ooo);(\psi_\vp(\na u)+\vp\na u)\cdot\nu=0\mbox{ on }\pp\ooo\},\earr\end{equation}where $\psi_\vp:\rr^d\to\rr^d$ is given by
\begin{equation}\label{e4.15a}
\psi_\vp(u)=\left\{\barr{ll}
\dd\frac1{\vp}\,v&\mbox{for }|v|\le\vp,\vsp
 \dd\frac v{|v|}&\mbox{for }|v|>\vp.\earr\right.\end{equation} It is easily seen that $\pp\vf_\vp$  (respectively, $\vf_\vp)$ is an approximation of $\pp\vf$ (res\-pec\-ti\-vely, $\vf$ in sense of graph convergence.

 In fact, we have

 \begin{proposition}\label{p4.1} For each $h>0$ and $f\in L^2(\ooo)$, we have
\begin{equation}\label{e4.16}
\lim_{\lbb\to0}(I+h\pp\vf_\vp)^{-1}f=(I+h\pp\vf)^{-1}f\end{equation}in the strong topology of $L^2(\ooo)$.\end{proposition}

\n{\bf Proof.} We set $u_\vp=(I+h\pp\vf_\vp)^{-1}f$. We have, by \eqref{e4.15},
\begin{equation}\label{e4.17}
\barr{ll}
u_\vp-h\vp\D u_\vp-h\,{\rm div}\,\psi_\vp(\na u_\vp)=f&\inr\ \ooo,\vsp
u_\vp\cdot\nu=0&\onr\ \pp\ooo.\earr\end{equation}By \eqref{e4.15a}, \eqref{e4.17}, we get the \fwg\ estimates (for fixed $h>0$)
\begin{eqnarray}
\int_\ooo\(\vp|\na u_\vp|^2+\frac12\,u^2_\vp+j_\vp(|\na u_\vp|)\)dx&\le&\frac12\int_\ooo f^2dx\label{e4.18}\\[1,2mm]
|\psi_\vp(\na u_\vp)|_\9&\le&1,\ \ \ff\vp>0\label{e4.19}\\[1,2mm]
\int_\ooo|\na u_\vp|dx&\le& C,\ \ \ff\vp>0.\label{e4.20}
\end{eqnarray}We also have
\begin{equation}\label{e4.21b}
-\int_\ooo \D u_\vp\,{\rm div}\,\psi_\vp(\na u_\vp)dx\ge0,\ \ \ff\vp>0,\end{equation}and this yields
\begin{equation}\label{e4.21}
\vp\int_\ooo|\D u_\vp|^2dx\le C,\ \ \ff\vp>0.\end{equation}We note that \eqref{e4.21b} is a direct consequence of the inequality
$$\int_\ooo|\na(I-\vp\D)^{-1}u|dx\le\int_\ooo|\na u|dx,$$for all $u\in H^2(\ooo);\ \frac{\pp u}{\pp n}=0$ on $\pp\ooo$. The last inequality essentially due to H.~Brezis is established in \cite{8a} (see Proposition 8.1) for $u\in H^2(\ooo)\cap H^1_0(\ooo)$ but it remains true in the present case too.

Moreover, taking into account that
$$\barr{l}
(\psi_\vp(\na u_\vp)-\psi_\mu(\na u_\mu))(\na u_\vp-\na u_\mu)\ge\vp|\psi_\vp(\na u_\vp)|^2+\mu|\psi_\mu(\na u_\mu)|^2\earr$$we see by \eqref{e4.17}--\eqref{e4.20} that $\{u_\vp\}$ is a Cauchy sequence in $L^2(\ooo)$. Hence, for $\vp\to0$,
$$\barr{rcll}
u_\vp&\longrightarrow&u&\mbox{strongly in $L^2(\ooo)$},\vsp
\psi_\vp(\na u_\vp)&\longrightarrow&\eta&\mbox{weak-star in $L^\9(\ooo)$},\vsp
{\rm div}\,\psi_\vp(\na u_\vp)&\longrightarrow&{\rm div}\,\eta&\mbox{weakly in $L^2(\ooo)$}.\earr$$We have, therefore,
\begin{equation}\label{e4.21a}
u-h\,{\rm div}\,\eta=f\ \inr\ \ooo.\end{equation}Moreover, by the inequality
$$\barr{l}
\dd\int_\ooo(j_\vp(\na u_\vp)-j_\vp(\na v))dx\le\dd\int_\ooo\psi_\vp(\na u_\vp)(\na u_\vp-\na v)dx\vsp
\qquad=-\dd\int_\ooo{\rm div}\psi_\vp(\na u_\lbb)(u_\vp-v)dx,\ \ff v\in H^1(\ooo),\earr$$we obtain that
\begin{equation}\label{e4.22}
-\int_\ooo{\rm div}\,\eta(u-v)dx\ge\vf(u)-\vf(v),\ \ff v\in BV(\ooo),\end{equation}because, by \eqref{e4.20}, it follows that (see Proposition 10.1.1 in \cite{5a})
$$\int_\ooo|Du|\le\liminf_{\vp\to0}\int_\ooo|\na u_\vp|dx$$and$$\left|j_\vp(\na u)-|\na u|\,\right|\le\frac12\,\vp\mbox{\ \ a.e. in }\ooo.$$By \eqref{e4.21}, \eqref{e4.22}, we see that $u=(I+h\pp\vf)^{-1}f$, as claimed.

\begin{corollary}\label{c4.1} Let $u_0\in L^2(\ooo)$ and $u_\vp\in C([0,T];L^2(\ooo))\cap W^{1,\9}([\de,T];H)\cap L^2(\de,T;H^2(\ooo))$, $\ff \de\in(0,T)$, be the solution to the Cauchy problem
\begin{equation}\label{e4.23}\barr{ll}
(\psi_\vp(\na u_\vp)+\vp\na u_\vp)\cdot\nu=0&\onr\ (0,T)\times\pp\ooo,\vsp
u_\vp(0)=u_0&\inr\ \ooo,\vsp
\dd\frac{\pp u_\vp}{\pp t}-{\rm div}(\psi_\vp(\na u_\vp))-\vp\D u_\vp=0&\inr\ (0,T)\times\ooo.\earr\end{equation}
Then, for $\vp\to0$, $u_\vp(t)\to u(t)$ in $L^2(\ooo)$ is uniformly on $[0,T]$. Moreover, if $u_0\in H^1(\ooo)$, then $u\in L^\9(0,T;H^1(\ooo)),$ $\frac{du}{dt}\in L^\9(0,T;(H^1(\ooo))').$\end{corollary}

\n{\bf Proof.} By the Trotter-Kato theorem for nonlinear semigroups of contractions (see, e.g., \cite{8}, p.~168), it follows by Proposition \ref{p4.1} that
\begin{equation}\label{e4.25a}
\lim_{\vp\to0}e^{t\pp\vf_\vp}u_0=e^{t\pp\vf}u_0\mbox{ \ \ uniformly on $[0,T]$ in $L^2(\ooo)$},\end{equation}and this completes the proof of the first part. Assume now that $u_0\in H^1(\ooo)$. Then, by \eqref{e4.17} and \eqref{e4.21b}, it follows that
$$\|\na u_\vp\|_{L^2(\ooo)}\le\|\na f\|_{L^2(\ooo)},\ \ff f\in H^1(\ooo),$$and, therefore, by \eqref{e4.25a} we see that
$$\|u(t)\|_{H^1(\ooo)}\le\|u_0\|_{H^1(\ooo)},\ \ff t\in[0,T],$$as claimed.\mk

The constant $\lbb>0$ arising in the minimization problem \eqref{e4.7} is a "scale" parameter which in the corresponding dynamic model \eqref{e4.13} plays the role of a time discrete variable $\frac 1h$. As regards the $L^2$-square norm $\int_\ooo|u-u_0|^2dx$, it is a so called "fidelity" term, which has the origin in the least square approach  of inverse ill posed problems. Several authors proposed to replace this term by the $L^1$-norm of $u-u_0$ and so the energy functional \eqref{e4.2} (or, more exactly, \eqref{e4.4}) reduces to
\begin{equation}\label{e4.24}
{\rm Min}\left\{\int_\ooo|Du|+\lbb\int_\ooo|u-u_0|dx\right\}.\end{equation}The Euler--Lagrange optimality equations in this case are of the form
\begin{equation}\label{e4.25}
\pp\vf(u)+\lbb\,{\rm sgn}(u-u_0)\ni0,\end{equation}where ${\rm sgn}\,v=\frac v{|v|}$ for $v\ne0,\ {\rm sgn}\,0=[-1,1].$

As regards the existence in \eqref{e4.25}, we have

\begin{theorem}\label{t4.3} Let $1\le d\le2$. Then, for each $u_0\in L^1(\ooo)$, equation \eqref{e4.25} has at least one solution $u\in BV(\ooo)$.\end{theorem}

\n{\bf Proof.} The operators $\pp\vf$ and $B(u)=\lbb\,{\rm sgn}(u-u_0)$ are maximal monotone in $L^2(\ooo)$ and, since $D(B)=L^2(\ooo)$, so is $B+\pp\vf$ (see, e.g., \cite{8}, p.~43).  Hence, for each $\vp>0$, the equation
$$\vp u_\vp+\pp\vf(u_\vp)+\lbb\,{\rm sgn}(u_\vp-u_0)\ni0$$has a unique solution $u_\vp\in D(\pp\vf)$. We have the estimates
$$\vp|u_\vp|^2_2+\int_\ooo|Du_\vp|\le C,\ \ \ff \vp>0,$$and  the embedding $BV(\ooo)\subset L^2(\ooo)$. It follows, therefore, that on a subsequence, again denoted $\vp$, we have
$$\barr{lcll}
 \vp u_\vp&\to&0&\mbox{\ \ strongly in $L^2(\ooo)$}\vsp
 u_\vp&\to& u&\mbox{\ \ weakly in $L^2(\ooo)$}.\earr$$Since $\pp\vf+B$ is maximal monotone in $L^2(\ooo)\times L^2(\ooo)$, we infer that
$$\pp\vf(u)+\lbb\,{\rm sgn}(u-u_0)\ni0,$$as claimed.

A solution $u_\9$ to   \eqref{e4.25} can be obtained by letting $t\to\9$ in the evolution equation
$$\barr{l}
\dd\frac{du(t)}{dt}+\pp\vf(u(t))+\lbb\,{\rm sgn}(u(t)-u_0)\ni0,\ t\ge0,\vsp
u(0)=0.\earr$$Equivalently,
\begin{equation}\label{e4.26}
\barr{l}
\dd\frac{du}{dt}+\pp\Phi(u)\ni0,\ \ t\ge0,\vsp
u(0)=0,\earr\end{equation}where $\Phi(u)=\vf(u)+\lbb|u-u_0|_1.$

We know that (see \cite{8}, p. 172), for $t\to\9$,
\begin{equation}\label{e4.27}
u(t)\to u_\9\mbox{\ \ weakly in $L^2(\ooo)$},\end{equation}where $u_\9$ is a solution to \eqref{e4.25}. On the other hand, as easily seen by \eqref{e4.26}, we have
$$\Phi(u(t))+\int^t_0\left|\frac{du}{ds}\,(s)\right|^2_2ds\le\Phi(0).$$
This means that, if $d=1$, then the orbit $\{u(t);\ t\ge0\}$ is compact in $L^2(\ooo)$ and so \eqref{e4.27} implies
\begin{equation}\label{e4.28}
\lim_{t\to\9}u(t)=u_\9\mbox{\ \ strongly in }L^2(\ooo).\end{equation}In particular, this implies that the steepest descent scheme
$$u_{i+1}=u_i-h\pp \Phi(u_{i+1}),\ i=0,1,...,$$is strongly convergent to a solution $u_\9$ to \eqref{e4.25}.

 To summarize, one can say that, due to the nonlinear semigroup theory, one has a rigorous mathematical theory for the total variation model in image processing and, from this perspective, the contributions given in the works \cite{3}, \cite{4}, \cite{5} should be also  emphasized.

 \begin{remark}\label{r4.1} {\rm As seen in Corollary \ref{c4.1}, if $u_0\in H^1_0(\ooo)$ (and this extends to $u_0\in H^1(\ooo)$), then $u(t)\in H^1_0(\ooo)$ a.e., $t>0$. However, for denoising purposes condition $u_0\in H^1_0(\ooo)$ is not realistic because the blurred image $u_0$ is less regular.}\end{remark}

\section{The Perona-Malik model revisited}
\setcounter{theorem}{0}
\setcounter{equation}{0}

We come back to the Perona--Malik equation  \eqref{e1.11}--\eqref{e1.12}. As noticed \mbox{earlier,} this equation is ill posed and  in literature there are several attempts to circumvent this inconvenience. Typically in most situations one re\-gu\-la\-rizes the diffusivity $u\to g(|\na u|^2) $ as in \eqref{e1.13},  for instance. This procedure transforms \eqref{e1.11} in a well posed problem, but affects the filtering property and, in special, the edge identification, because the regularized equation is well posed in a high order Sobolev space (in $H^2(\ooo)$. for instance).

For simplicity, we consider here problem \eqref{e1.11} with Dirichlet boundary conditions, that is,
\begin{equation}\label{e5.1}
\barr{ll}
\dd\frac{\pp u}{\pp t}-{\rm div}\,g((|\na u|^2)\na u)=0&\inr\ (0,T)\times\ooo,\vsp
u=0&\onr\ (0,T)\times\pp\ooo,\vsp
u(0)=u_0&\inr\ \ooo,\earr\end{equation}where $g$ is given by \eqref{e1.12}.

We have

\begin{theorem}\label{t5.1} Assume that $\ooo$ is a bounded convex set with smooth boundary. Let $u_0\in W^{1,\9}(\ooo)$. Then, for $\a>\|u_0\|_{W^{1,\9}(\ooo)}$, there is a unique solution $u$ to \eqref{e5.1} satisfying
\begin{eqnarray}
& u\in C([0,T];L^2(\ooo))\cap L^\9(0,T;H^1_0(\ooo)),\label{e5.2}\\[1,2mm]
 &\frac{\pp u}{\pp t},{\rm div}(g(|\na u|^2)\na  u)\in L^\9(\de,T;L^2(\ooo)),\ \ff\de\in(0,T),\label{e5.3}\\[1,2mm]
& |\na_xu(t,x)|\le\a,\ \mbox{a.e. }(t,x)\in(0,T)\times\ooo,\label{e5.3a}\\[1,2mm]
 &\dd\frac{\pp^+}{\pp t}\,u(t,x)-{\rm div}(g(|\na u(t,x)|^2)\na u(t,x))=0,\ \ff t\in(0,T),\ x\in\ooo.\ \ \label{e5.4}
\end{eqnarray}\end{theorem}

Taking  into account that $u_0$ is the blurred image, the condition $u_0\in W^{1,\9}(\ooo)$ is, apparently, too restrictive. However, the above result shows that the Perona-Malik model is well posed in the class of "moderately" degraded original images $u_0$. On the other hand, the smoothing effect of the model is underlined by property \eqref{e5.2}, which implies that $u(t,\cdot)\in H^1(\ooo)$ for all $t\in(0,T).$ However, the diffusion term in \eqref{e5.3} is attenuated  of order $|\na u(t)|^2$, which emphasizes  the edges detection performance of the model.

\bk\n{\bf Proof of Theorem \ref{t5.1}.} We set
$$g_\a(s)=\left\{\barr{ll}
g(s)=\dd\frac{\a^2}{s+\a^2}&\mbox{ for }0\le s\le\a,\vsp
\dd\frac\a{\a+1}&\mbox{ for }s>\a,\earr\right.$$and consider the operator $A^\vp_\a:D(A^\vp_\a)\subset L^2(\ooo)\to L^2(\ooo),$ $\vp>0$, defined by
\begin{equation}\label{e5.5}
A^\vp_\a u=-{\rm div}(g_\a(|\na u|^2)\na u)-\vp\D u,\ u\in D(A^\vp_\a),\end{equation}where $D(A^\vp_\a)=H^1_0(\ooo)\cap H^2(\ooo).$ We have

\begin{lemma}\label{l5.1} For each $\vp>0$, the operator $A^\vp_\a$ is maximal monotone\break $(m$-accre\-tive$)$ in $L^2(\ooo)\times L^2(\ooo)$.\end{lemma}

\n{\bf Proof.} Consider the operator $\wt A^\vp_\a:H^1_0(\ooo)\to H^{-1}(\ooo)$ defined by
$$_{H^{-1}(\ooo)}\<\wt A^\vp_\a u,v\>_{H^1_0(\ooo)}=\int_\ooo g_\a(|\na u|^2)\na u\cdot\na v\,dx+\vp\int_\ooo\na u\cdot\na v\,dx,\ \ff v\in H^1_0(\ooo).$$It is easily seen the the mapping
$r\to g_\a(|r|^2)r$ is monotone and continuous in $\rr^d$. This implies that the operator $\wt A^\vp_\a$ is monotone, demicontinuous (that is, strongly-weakly continuous) and coercive. Then, according to a well known result of Minty and Browder (see, e.g., \cite{8}, p~80), it is maximal monotone and surjective. Hence, its restriction $A^\vp_\a$ to $L^2(\ooo)$ is maximal monotone, that is, $\rr(I+A_\a^\vp)=L^2(\ooo)$. This means that, for each $u_0\in L^2(\ooo)$, the Cauchy problem
\begin{equation}\label{e5.6}
\barr{l}
\dd\frac{du}{dt}+A^\vp_\a u=0,\ \ t\in(0,T),\vsp
u(0)=u_0,\earr\end{equation}has a generalized solution $u_\vp\in C([0,T];L^2(\ooo))$ (see Theorem \ref{t2.1}).

Moreover, since, as easily seen, $A^\vp_\a=\na \vf^\vp_\a$, where
$$\barr{lcl}
\vf^\vp_\a(u)&=&\dd\int_\ooo(j_\a(|\na u|)+\vp|\na u|^2)dx,\vspp
j_\a(s)&=&\left\{\barr{ll}
\dd\frac{\a^2}2\,\log(s^2+\a^2)&\mbox{ for }0\le s\le\a,\vspp
\dd\frac{s^2}4+C(\a)&\mbox{ for }s>\a,\earr\right.\earr$$where $C(\a)=\frac{\a^2}4(\log\a+2\log2-1),$ it follows that $u_\vp$ is a strong solution of \eqref{e5.6} on $(0,T)$ (see Theorem \ref{t2.2}) and
\begin{equation}\label{e5.7}
\frac{du_\vp}{dt}\in L^\9(\de,T;L^2(\ooo)),\ u_\vp\in L^\9(\de,T;H^1_0(\ooo)\cap H^2(\ooo)),\ \ff \de\in(0,T).\end{equation}

Consider, now, the set closed convex $K\subset L^2(\ooo)$
$$K=\{u\in W^{1,\9}(\ooo)\cap H^1_0(\ooo);|\na u|_\9\le\a\}.$$We have

\begin{lemma}\label{l5.2} If $u_0\in K$, then $u_\vp(t)\in K$, $\ff t\in[0,T].$\end{lemma}

\n{\bf Proof.} According to a well known invariance theorem for the semigroups of nonlinear contractions (see H. Brezis \cite{9}), it suffices to show that, for each $\lbb>0$, the resolvent $(I+\lbb A^\vp_\a)^{-1}$ leaves invariant the set $K$. In other words, for each $f\in K$, the solution $v\in H^1_0(\ooo)$ to the equation
$$v-\lbb\,{\rm div}(g_\a(|\na v|^2)\na v)-\lbb\vp\D v =f\ \inr\ \ooo,$$belongs to $K$. Since $\ooo$ is convex, this follows by Theorem III.1 in \cite{10} and so Lemma \ref{l5.2} follows.

\bk\n{\bf Proof of Theorem \ref{t5.1} (continued).} We come back to equation\eqref{e5.6} and show that, if $u_0\in K$, then, for $\vp\to0$, $u_\vp$ is convergent to a solution $u$ to \eqref{e5.1}. Since $g_\a=g$ on $K$, we have
\begin{equation}\label{e5.8}
\barr{l}
\dd\frac{\pp u_\vp}{\pp t}-{\rm div}(g(|\na u_\vp|^2)\na u_\vp)-\vp\D u_\vp=0\ \inr\ (0,T)\times\ooo,\vsp
u_\vp(0)=u_0\ \inr\ \ooo,\ \ u_\vp=0\ \onr\ (0,T)\times\pp\ooo.\earr\end{equation}We note that,  taking into account that $g_\a(|\na u|^2)\na u=\pp j_\a(|\na u|)$, we have by Lemma  5.1  in \cite{8a} that
$$\int_\ooo{\rm div}(g_\a(|\na v|^2)\na v)\D v\,dx\ge0,\ \ff v\in H^1_0(\ooo)\cap H^2(\ooo).$$Then, by multiplying \eqref{e5.8} with $u_\vp,\D u_\vp$ and $\frac{\pp u_\vp}{\pp t}$, we get the apriori estimates
\begin{eqnarray}
\int^T_0\int_\ooo\left|\frac{\pp u_\vp}{\pp t}\right|^2dt\,dx+ \vf^\vp_\a(u_\vp(t)) &\le&\vf^\vp_\a(u_0)\le C,\label{e5.9}\\[2mm]
\frac12\,|u_\vp(t)|_2+\vp\int^T_0|\D  u_\vp(t)|^2_2dt&\le&\frac12|u_0|^2_2,\label{e5.10}
\end{eqnarray}
and, by Lemma \ref{l5.2},
\begin{eqnarray}
|\na u_\vp(t)|_\9&\le&\a,\mbox{\ \ a.e. }t\in(0,T).\label{e5.11}
\end{eqnarray}In particular, this implies that $\{u_\vp\}$ is compact in $C([0,T];L^2(\ooo))$ and, therefore, there is $u\in C([0,T];L^2(\ooo))\cap W^{1,2}([0,T];L^2(\ooo))$ such that, for $\vp\to0$,
$$\barr{rcll}
u_\vp&\longrightarrow&u&\mbox{strongly in }C([0,T];L^2(\ooo))\\
&&&\mbox{weakly in }L^2(0,T;H^1_0(\ooo)),\vspp
\na u_\vp&\longrightarrow&\na u&\mbox{weak-star in }L^\9((0,T)\times\ooo),\vspp
\vp\D u_\vp&\longrightarrow&0&\mbox{strongly in }L^2(0,T;L^2(\ooo)),\vspp
{\rm div}\,g(|\na u_\vp|^2)\na u_\vp)&\longrightarrow&\eta&\mbox{weakly in }L^2(0,T;L^2(\ooo)),\earr$$where
\begin{equation}\label{e5.12}
\barr{ll}
\dd\frac{\pp u}{\pp t}-\eta=0&\mbox{ a.e. }t\in(0,T), \ x\in\ooo,\vsp
u(0)=u_0&\inr\ \ooo ,\vsp
|\na u(t,x)|\le\a,&\mbox{ a.e. }(t,x)\in(0,T)\times\ooo.\earr\end{equation}To complete the proof, it suffices to show that
\begin{equation}\label{e5.13}
\eta={\rm div}(g(|\na u|^2)\na u)\ \inr\ \cald'((0,T)\times\ooo).\end{equation}
$(\cald'((0,T)\times\ooo)$ is the space of distributions on
$(0,T)\times\ooo.)$

To this end, we recall that $A^\vp_\a=\na\vf^\vp_\a$ and this implies that, for all $v\in L^2(0,T;H^{-1}_0(\ooo))$,
$$\dd\int^T_0\int_\ooo{\rm div}(g(|\na u_\vp|^2)\na u_\vp)(u_\vp-v)dx\,dt\ge
\dd\int^T_0\int_\ooo(j_\a(|\na u_\vp|)-j_\a(|\na v|))dx\,dt .$$
Letting $\vp\to0$, we obtain, for all $  v\in L^2(0,T;H^1_0(\ooo))$,
\begin{equation}\label{e5.14}
\int^T_0\int_\ooo\eta(u-v)dx\,dt\ge\int^T_0\int_\ooo(j_\a(|\na u|)-j_\a(|\na v|))dx\,dt,\end{equation} because the function $u\to\int^T_0\int_\ooo j_\a(|\na u|)dx\,dt$ is convex, lower semi\-con\-ti\-nuous and,  therefore, weakly lower semi\-con\-ti\-nuous in the space $L^2(0,T;H^1_0(\ooo)).$  By \eqref{e5.14}, it follows that $\eta\in\pp\vf_\a(u)$, where $\vf_\a:L^2(\ooo)\to\ov\rr$ is defined by
$$\vf_\a(u)=\int_\ooo j_\a(|\na u|)dx,\ u\in H^1_0(\ooo),$$and so \eqref{e5.13} follows. This completes the proof of the existence.

The uniqueness follows by the $L^2(\ooo)$ monotonicity of the operator $$u\buildrel A\over\longrightarrow-{\rm div}(g(|\na u|^2)\na u)$$ on the set $\{u\in H^1_0(\ooo);\ |\na u|_\9\le\a\}.$ Theorem \ref{t5.1} amounts to saying that the operator $u\buildrel A\over\longrightarrow-{\rm div}(g(|\na u|^2)\na u)$ with appropriate homogenous boundary conditions (Dirichlet or Neumann) generates a semiflow $t\to u(t,u_0)$ on the set $K$, which can be computed via the finite difference scheme
\begin{equation}\label{e5.15}
\barr{ll}
u_{i+1}-h\,{\rm div}(g(|\na u_{i+1}|^2)\na u_{i+1})=u_i& \inr\ \ooo,\ i=0,1,...,\vsp
u_{i+1}=0&\onr\ \pp\ooo, \earr \end{equation}or, equivalently, by the exponential formula

\begin{equation}\label{e5.16}
u(t,u_0)=\lim_{t\to\9}\(I+\frac tn\,A\)^{-n}u_0\mbox{ uniformly on }[0,T].\end{equation}

In order to extend the algorithm \eqref{e5.15} to all $u_0\in L^2(\ooo)$, that is, to general blurred images $u_0$, we consider the projection $P_K:L^2(\ooo)\to K$ and replace \eqref{e5.15} by
\begin{equation}\label{e5.17}
\barr{l}
u_{i+1}-h\,{\rm div}\,g(|\na u_{i+1}|^2\na u_{i+1})=u_i\ \inr\ \ooo,\ i=1,...,\vsp
u_1=P_Ku_0,\ \ u_{i+1}=0\ \onr\ \pp\ooo.\earr\end{equation}The projection $v=P_K f$, for $f\in L^2(\ooo)$, is given by
\begin{equation}\label{e5.18}
\barr{ll}
v-{\rm div}\,\b(\na v)\ni f&\inr\ \ooo,\vsp
v=0&\onr\ \pp\ooo,\earr\end{equation}where $\b:\rr^d\to 2^{\rr^d}$ is the normal cone to the ball $\{r\in\rr^d;\ |r|\le\a\}$, that is,
$$\b(r)=\left\{\barr{cl}
\lbb\,\frac r{|r|},&\mbox{ for }|r|=\a,\vsp
0&\mbox{ for }|r|<\a,\earr\right.$$where $\lbb>0$.

The equation \eqref{e5.18}, which is well posed in $H^1_0(\ooo)$ and is equivalent with
$${\rm Min}\left\{\int_\ooo(v-f)^2dx;\ |\na v(x)|\le\a,\mbox{ a.e. }x\in\ooo\right\}$$can be approximated by the generalized problem
$$v_\vp-{\rm div}\,\b_\vp(\na v_\vp)=f\ \inr\ \ooo,\ v_\vp=0\ \onr\ \pp\ooo,$$where
$$\b_\vp(v)=\frac1\vp\,\frac v{|v|}\,(|v|-\a)^+,\ \ff v\in\rr^d.$$Then, \eqref{e5.17} leads to the \fwg\ denoising algorithm
\begin{eqnarray}
&&\barr{l}
u^\vp_{i+1}-h\,{\rm div}(g(|\na u^\vp_{i+1}|^2)\na u^\vp_{i+1})=u^\vp_i\ \inr\ \ooo,\ i=1,2,...\vsp
u^\vp_{i+1}=0\ \onr\ \pp\ooo,\earr\label{e5.19}\\[2mm]
&&\barr{l}
u^\vp_1-{\rm div}\,\b_\vp(\na u^\vp_1)=u_0\ \inr\ \ooo,\ u_1^\vp=0\ \onr\ \pp\ooo.\earr\label{e5.20}
\end{eqnarray}

\section{A denoising model based on the porous\\ media equation}
\setcounter{theorem}{0}
\setcounter{equation}{0}

The method we discuss here and developed in \cite{7}  starts from a variational problem considered in the Sobolev distribution space $H^{-1}(\ooo)$ for which the corresponding Euler--Lagrange equation is a nonlinear elliptic diffusion equation related to porous media equations. An important feature of this approach is that the observation $f$ can be taken as a distribution with support in a finite number of points $\{(x_1)_\ell,(x_2)_h)\}_{\ell,h}\subset\rr^2$. This means that it can be applied to the inpainting problem for a subdomain $\wt\ooo$ of $\ooo$ with zero Lebesgue measure and, in particular, to the case where $\wt\ooo=\bigcup^N_{j=1}\G_j$, where $\G_j$ are smooth curves.

Denote by $H^{-1}(\ooo)$ the dual space of $H^1_0(\ooo),$ that is, the space of all distributions $u\in\cald'(\ooo)$ of the form
$$u=\sum^d_{j=1}\frac{\pp w_j}{\pp x_j}\ \inr\ \cald'(\ooo),\ w_j\in L^2(\ooo).$$The norm of $H^{-1}(\ooo)$ is denoted by $\|\cdot\|_{-1}$ and we have
$$\|u\|^2_{-1}=\int_\ooo(-\D)^{-1}u(x)u(x)dx,\ \ff u\in H^{-1}(\ooo),$$where $(-\D)^{-1}u=w$ is the solution to the Dirichlet problem
$$-\D w=u\ \inr\ \ooo;\ \ w=0\ \onr\ \pp\ooo.$$The scalar product in $H^{-1}(\ooo)$ is given by $$\<u,v\>_{-1}=\int_\ooo(-\D)^{-1}u(x)v(x)dx.$$

As in the previous situations, we identify an image with a real valued function $u$ defined on an open bounded domain $\ooo\subset\rr^2$. At each point $x\in\ooo$, $u(x)$ is the gray  level of image at $x$. Usually, $u(x)$ are coded with integer values in the interval $[0,255]$, where $0$ represents the {\it black level} and $255$ the {\it white level}.

The problem to be considered here is to restore the image $u$ from the observation$f_j$ of the blurred image $\wt\ooo$. Namely, $f\in H^{-1}(\ooo)$ might be defined~by
$$f(\vf)=\sum^M_{j=1}\int_{\G_j}f_j(s)\vf(s)ds,\ \ \ff\vf\in H^1_0(\ooo),$$where $f_j\in L^2(\G_j).$

The restored image $u$ is found from the minimization problem
\begin{equation}\label{e6.1}
{\rm Min}\!\left\{\int_\ooo \! j(u(x))dx
+\frac12\,\|u{-}f\|^2_{-1};
\, u\in L^1(\ooo),\, u{-}f\in H^{-1}(\ooo)\right\},\end{equation}where $j:\rr\to\rr$ is a convex and lower semicontinuous function.

We denote by $\b:\rr\to2^\rr$ the subdifferential of $j$, that is,
$$\b(r)=\{\oo\in\rr;\ w(r-s)\ge j(r)-j(s),\ \ff s\in\rr\}.$$

Formally, the Euler--Lagrange optimality conditions in \eqref{e6.1} are given by the elliptic boundary value problem
\begin{equation}\label{e6.2}
-\D\b(u)+u=f\ \inr\ \ooo,\ \ \ \b(u)=0\ \onr\ \pp\ooo,\end{equation}but the exact meaning of \eqref{e6.2} is made precise below.

\begin{theorem}\label{t6.1} Assume that $f\in H^{-1}(\ooo)$ and that $j$ is convex, lower semicontinuous such that $j(0)=0$ and
\begin{equation}\label{e6.3}
\lim_{|r|\to\9}\frac{j(r)}{|r|}=+\9.\end{equation}Then, the minimization problem \eqref{e6.1} has a unique solution $u^*$ which sa\-tis\-fies equation \eqref{e6.2} in the \fwg\ sense. There is $\eta^*\in H^1_0(\ooo)$ such that $\eta^*(x)\in\b(u^*(x)),$ a.e. $(x,y)\in\ooo$ and
\begin{equation}\label{e6.4}
\int_\ooo\na\eta^*\cdot\na\vf\,dx
+{}_{H^{-1}(\ooo)}\!\!\<u^*-f,\vf\>_{H^1_0(\ooo)}=0,\ \ff\vf\in H^1_0(\ooo).
\end{equation}
\end{theorem}

\n{\bf Proof.} The functional
$$\Phi(u)=\int_\ooo j(u)dx\,dy+\frac12\,\|u-f\|^2_{-1}$$is convex and lower semicontinuous on $H^{-1}(\ooo)$ and its subdifferential $\pp\Phi: H^{-1}(\ooo)\to H^{-1}(\ooo)$ is given by
$$\pp\Phi(u)=\pp\psi(u)+u-f,$$where $\pp\psi$ is the subdifferential of the function $\psi:H^{-1}(\ooo)\Longrightarrow\ov\rr$
$$\psi(u)=\int_\ooo j(u)dx,\ u\in L^1(\ooo)\cap H^{-1}(\ooo).$$On the other hand, according to a well known result due to Brezis (see \cite{8}, p.~76), the subdifferential $\pp\psi:H^{-1}(\ooo)\to H^{-1}(\ooo)$ is given by
$$\barr{r}
\dd\pp\psi(y)=\{w\in H^{-1}(\ooo);\ w=-\D\eta,\ \eta\in H^1_0(\ooo),\ \eta(x,y)\in\b(u(x,y))\\\mbox{ a.e. }(x,y)\in\ooo\}.\earr$$Hence,
$$\pp\Phi(u)=\{-\D\eta+u-f;\ \eta\in\b(u)\mbox{ a.e. in }\ooo\}$$and, therefore, the minimization problem \eqref{e6.1} reduces to $\pp\Phi(u)=0$, that is, to \eqref{e6.4}.

On the other hand, since the function $\Phi$ is convex, lower semicontinuous and coercive on $H^{-1}(\ooo)$, there is $u^*$ such that
$$u^*={\rm arg}\,\min\{\Phi(u);\ u\in H^{-1}(\ooo)\}.$$Clearly, $u^*$ is unique and belongs to $L^1(\ooo)$. This completes the proof.

\begin{remark}\label{r2.1} {\rm We see in Theorem \ref{t6.2} below that Theorem \ref{t6.1} extends to $f\in M(\ov\ooo)$, the space of bounded regular measures on $\ov\ooo_0$. }\end{remark}

If condition \eqref{e6.3} does not hold, then the functional $\Phi$ is no longer lower semicontinuous on $H^{-1}(\ooo)$ and so the minimum point in \eqref{e6.1} no longer exists. However, we still have in this case existence of a weak solution to the elliptic boundary value problem \eqref{e6.4}.

\begin{theorem}\label{t6.2} Assume that $f\in M(\ov\ooo)$ and that $j:\rr\to]-\9,+\9]$ is convex and lower semicontinuous, then there is a unique pair of functions$$(u^*,\eta^*)\in L^1(\ooo)\times W^{1,p}_0(\ooo),\ 1\le p<2,$$which satisfies \eqref{e6.4} and such that
$$u^*-f\in W^{-1,p}(\ooo),\ \ \frac1{p'}+\frac1p=1.$$\end{theorem}
Equation \eqref{e6.4} is considered here in the \fwg\ weak sense
\begin{equation}\label{e6.5}
\int_\ooo\na \eta^*\cdot\na\vf\,dx +{}_{W^{-1,p}(\ooo)}\!\<u-f,\vf\>_{W^{1,p'}_0(\ooo)}=0,\ \ff\vf\in W^{1,p'}_0(\ooo),\end{equation}where $W^{-1,p}(\ooo)$ is the dual of the Sobolev space $W^{1,p'}_0(\ooo).$ We note that Theorem \ref{t6.2} is, in particular,  applicable to the case where $f$ is of the form
$$f=\sum^M_{j=1} \mu_{j}\de(x_j)$$and $\de(x_j)$ is the Dirac measure concentrated in $ x_j $. This corresponds to the case where the image is observable in a finite set of pixels $\{ x_j \}.$

\bk\n{\bf Proof of Theorem \ref{t6.2}.} We set $\g=\b^{-1}$ and rewrite \eqref{e6.2} as
\begin{equation}\label{e6.6}
-\D v+\g(v)\ni f\ \inr\ \ooo,\ \ \ v=0\ \onr\ \pp\ooo.\end{equation}Next, we consider a sequence $\{f_j\}\subset L^2(\ooo)$ such that $f_j\to f$ in $M(\ov\ooo)$. For each $f=f_j$, by standard elliptic existence theory, \eqref{e6.6} has a unique solution $v_j\in H^1_0(\ooo)\cap H^2(\ooo)$ and the \fwg\ estimates hold
$$\|f_j\|_{L^1(\ooo)}\le C,\ \ff j,\ \ \ \|\g(v_j)\|_{L^1(\ooo)}\le\|f_j\|_{L^1(\ooo)}\le C,\ \ff j.$$Therefore, $\{\D v_j\}$ is bounded in $L^1(\ooo)$. This implies that (see, e.g., \cite{8} or
\cite{9a}, p.~108) the sequence $\{v_j\}$ is bounded in $W^{1,p}_0(\ooo)$ for $1\le p<2$, and so it is compact in $L^1(\ooo)$. Moreover, by the equation
\begin{equation}\label{e6.6a}
-\D v_j+u_j=f_j\ \inr\ \ooo;\ \ v_j=0\ \onr\ \pp\ooo,\end{equation}where $u_j\in \g(v_j)$, a.e. in $\ooo$, we have, by multiplying with $\zeta_{j,k}\in{\rm sgn}(u_j-u_k)$ and integrating on $\ooo$,
$$\|u_j-u_k\|_{L^1(\ooo)}
\le\|f_j-f_k\|_{L^1(\ooo)}\le\|f_j-f_k\|_{M(\ov\ooo)}.$$
Hence, there are $u^*,\eta^*$ such that
$$\barr{lcll}
u_f&\longrightarrow&u^*&\mbox{ strongly in $L^1(\ooo)$ and }\vsp
v_j&\longrightarrow&\eta^*&\mbox{ weakly in $W^{1,p}_0(\ooo)$ and strongly in $L^1(\ooo)$}\earr$$as $j\to\9$. We also have that

\begin{equation}\label{e6.7}
\int_\ooo\na v_j\cdot\na\vf\,dx\,dy+\int_\ooo(u_j-f_j)\vf\,dx\,dy=0,\ \ff\vf\in W^{1,p'}_0(\ooo).\end{equation}Since $W^{1,p'}_0(\ooo)\subset C(\ov\ooo)$ for $\frac1{p'}=1-\frac1p$, we may pass to limit into \eqref{e6.6a} to conclude that $(u^*,\eta^*)$ satisfies \eqref{e6.5}. The uniqueness is immediate by the monotonicity of $\b$. Thus, the proof is complete.\bk

We note that Theorem \ref{t6.2} is, in particular, applicable to $j(r)=|r|,$ $\b(r)={\rm sign}\,r$. However, in applications it is more convenient to take $\b$ a monotonically increasing $C^1$-function. Indeed, in this case,  the solution $u^*$ to \eqref{e6.2} satisfies
$$\na\b(u^*)=\b'(u^*)\na u^*\in L^p(\ooo),\ 1\le p<2,$$and this implies that $\b(u^*)$ is smooth (in fact, absolutely continuous on each horizontal or vertical line) but the gradient $|\na u^*|$ might become big in the region where $\b'(u^*)$ (or $u^*)$ is small.
A formal conclusion from this is that the above procedure has a smoothing effect on the initial image which might be extremely blurred, and it is also efficient to edge detection (that is, of sharp contrast). Also, we may take $\b$ of the form $\b(x,u)$ which leads to a denoising approach is that the solution $u$ to \eqref{e6.1} has an anisotropic diffusion. Of course, the denoising procedure is more  effective if $\b$ is closer of sign or $j(r)$ of $|r|$. However, in this limit case, the restored image looses its smoothness and so perhaps one of the bet choice might be $\b(r)=|r|^ar$, where $0<\a<1$, which is just the case covered by Theorem \ref{t2.1}.

As mentioned earlier, this model allows observation $f$ which are not standard $L^1$-functions on $\ooo$, but distributions (measures or $H^{-1}$-distributions). This is applicable, in particular, in the case of restoration of an image known on a set of strips with arbitrarily small width or on a finite set of points (pixels). (See \cite{7} for restoring the experiment in such a case.)

The dynamic model corresponding to \eqref{e6.1} is
\begin{equation}\label{e6.8}
\barr{ll}
\dd\frac{\pp u}{\pp t}-\D\b(u)\ni0&\inr\ (0,\9)\times\ooo,\vsp
u(0,x)=f(x),&x\in\ooo,\vsp
u=0&\onr\ (0,\9)\times\pp\ooo.\earr\end{equation}By Theorem \ref{t2.2}, it follows that, for each $u_0\in\ov{D(\pp\vf)}$, that is, $u_0\in L^1(\ooo)\cap H^{-1}(\ooo)$, the Cauchy problem \eqref{e6.8} has a unique solution $u\in C([0,T];H^{-1}(\ooo))$ with
$$\frac{du}{dt},\ \D\b(u)\in L^\9([\de,T];H^{-1}(\ooo)),\ \ff\de\in [0,T].$$
Moreover, the finite difference scheme
\begin{equation}\label{e6.9}
\barr{l}
u_{i+1}-h\D\b(u_{i+1}=u_i,\ \inr\ \ooo,\ i=0,...,\vsp
\b(u_{i+1})\in H^1_0(\ooo),\earr\end{equation}is convergent to the solution $u$.

By \eqref{e6.9}, we see that $\b(u_{i+1}\in H^1_0(\ooo)$ which emphasizes the smoothing effect of the model. On the other hand, the initial blurred image $u_0$ can be chosen in $L^1(\ooo)\cap H^{-1}(\ooo)$ which is the case with an extremely degraded image.
For related works, we cite also \cite{13}, \cite{14a}.

\bk\n{\bf Acknowledgement.} This work was supported by a grant of the Romanian National Authority for Scientific Research, CNCS-UEFISCDI Project PN-II-ID-PCE-2011-30--27.


\begin{thebibliography}{nn}

\bibitem{1} L.Alvarez, P.L. Lions, J.M. Morel, Image selective smoothing and edge detection by nonlinear diffusion, {\it SIAM J. Numer. Anal.}, 29 (1992), 845-866.

\bibitem{2}  L. Ambrosio, N. Fusco, D. Pallara, {\it Functions of Bounded Variations and Free Discontinuous Problems}, Oxford Mathematical Monographs, 2000.

\bibitem{3}  F. Andreu, C. Ballester, V. Caselles, J.M. Mazon, Minimizing total variation flow, {\it Differential and Integral Equations}, 14 (2001), 321-360.

\bibitem{4}  F. Andreu, V. Caselles, J.I. Diaz, J.M. Mazon, Qualitative properties of the total variation flow, {\it  J. Funct. Anal.}, 188 (2002), 516-547.

\bibitem{5}  F. Andreu,   J.M. Mazon, The total Variation Flow (Preprint).

\bibitem{5a} H. Attouch, G. Buttazzo, G. Michaille, {\it Variational Analysis in Sobolev and BV Spaces}, MPS-SIAM Series, Philadelphia, 2006.

\bibitem{8}  V. Barbu,  {\it Nonlinear Differential Equations of Monotone Type in Banach Spaces}, Springer, New York, Dordrecht, London, 2009.


\bibitem{7}  T. Barbu, V. Barbu, A PDE approach to image restoration problem with observation on a meager domain, {\it Nonlinear Anal.  Real World Appl.}, 13 (2012), 1206-1215.

\bibitem{6}  T. Barbu, V. Barbu, V. Biga, D. Coca, A PDE approach to image denoising and restoration, {\it Nonlinear Anal. Real. Word Appl.}, 10 (2009), 1351-1361.


    \bibitem{8a} V. Barbu, M. R\"ockner, Stochastic variational inequalities and applications to the total variation flow perturbed by linear multiplicative noise, {\it Archive Rational Mechanics and Analysis} (to appear).


        \bibitem{8aa} M. Black, C. Sapiro, D. Marimont, D. Heeger, Robust anisotropic diffusion, {\it IEEE Transactions on Image Procesing}, vol. 7, no. 3 (1998),   421-432.

   \bibitem{9a}  H. Brezis, {\it Functional Analysis, Sobolev Spaces and Partial Differential Equations}, Springer, New York, Dordrecht, London, 2010.

\bibitem{9}  H. Brezis, {\it Op\'erateurs maximaux monotones et semigroupes de contractions dans les espaces de Hilbert}, North-Holland, Amsterdam, 1973.

\bibitem{10} H. Brezis, G. Stampacchia, Sur la r\'egularit\'e  de la solution d'in\'equations elliptiques, {\it Bull. de la S.M.F.}, 96 (1968), 153-180.

\bibitem{10a} J. Calder, A. Mansouri, A. Yezzi, Image diffusion and sharpening via Sobolev gradient flows, {\it SIAM J. Imaging Sci.},  3 (2010), 981-1014.

    \bibitem{10aa} J. Calder, A. Mansouri, A. Yezzi, New possibilites in image diffusion and sharpening via high-order Sobolev gradients, {\it J. Math. Imaging},  40 (2011), 248-252.

\bibitem{11}  F. Catte, P.L. Lions, J.M. Morel, T. Col, Image selective smoothing and edge detection by nonlinear diffusion, {\it  SIAM J. Numerical Anal.}, 29 (1992), 182-193.

\bibitem{12}  A. Chamballe, P.L. Lions, Image recovery via total variation minimization and related problems, {\it Numer. Math.}, 76 (1997), 167-182.


    \bibitem{10aaa} T. Chan, J. Sheng, Mathematical models of local non-texture im\-pain\-ting, {\it SIAM J. Appl. Math.}, 62 (2001), 1019-1043.

    \bibitem{13aa} Z. Jin, X. Yang, Strong solutions for the generalized Perona--Malik equations for image restoration, {\it Nonlinear Anal.}, 73 (2010), 1077-1084.

    \bibitem{13a} J. Ka\v cur, K. Mikula, Slow and fast diffusion effects in image processing, {\it Comput. Visual. Sci.}, 3 (1996), 185-195.

\bibitem{13} J. Ka\v cur, K. Mikula, Solution of nonlinear diffusion appearing in image smoothing and edge detection, {\it Applied Numerical Mathematics}, 17 (1995), 47-69.

\bibitem{14} P.L. Lions, Axiomatic  derivation of image processing models, {\it Math. Models in Applied Sciences}, 4 (1994), 467-475.


    \bibitem{14a} S. Osher, A. Sole, L. Vese, Image decomposition and restoration \mbox{using} total variation minimization and the $H^{-1}$-norm, {\it Multiscale Model Si\-mu\-la\-tion}, 1 (2003), 349-370.

\bibitem{15} P. Perona, J. Malik, Scale space and edge detection using anisotropic diffusion, {\it IEEE Transactions Pattern Anal. Mach. Intell.}, 12 (1990), 629-639.

\bibitem{16} L. Rudin, S. Osher, E. Fatemi, Nonlinear total variation based noise removal algorithms, {\it Physica}, D. 60 (1992), 259-268.


\bibitem{18} J. Weickert, {\it Anisotropic Diffusion in Image Processing}, B.G. Teubner, Stuttgart, 1998.


\bibitem{17} J. Weickert, M. Bart, R. Haar, Efficient and reliable schemes for non\-linear diffusion filters, {\it IEEE Transactions on Image Processing}, 7 (1998), 398-410.


\end{thebibliography}
\end{document}